\documentclass[12pt]{article}
\usepackage{amsmath, amssymb, amsthm, enumerate}
\newtheorem{theorem}{Theorem}[section]
\newtheorem{lemma}[theorem]{Lemma}
\newtheorem{proposition}[theorem]{Proposition}
\newtheorem{remark}[theorem]{Remark}

\newtheorem{cor}[theorem]{Corollary}

\newtheorem{defin}[theorem]{Definition}

\def\rr{{\mathbb R}}

\def\cc{{\mathbb C}}

\def\qq{{\mathbb Q}}

\def\su{\subset}

\def\al{\alpha}
\def\be{\beta}

\def\De{\Delta}

\def\ph{\phi}

\def\cd{\cdot}
\def\stb{,\ldots ,}

\def\del{\partial}

\def\bsk{\bigskip}
\def\noi{\noindent}

\def\sumj0m{\sum_{j=0}^m}

\def\sumi0n{\sum_{i=0}^n}
\def\sumina{\sum_{i=1}^n a_i}

\def\suminf{\sum_{i=1}^n a_if(\al_ix+\be_iy)}
\def\suminftilde{\sum_{i=1}^n a_i\tilde{f}(\al_ix+\be_iy)}

\begin{tiny}\end{tiny}

\def\proof{\noi {\bf Proof.} }

\def\x1n{x_1 \stb x_n}
\def\y1n{y_1 \stb y_n}

\def\diag{{\rm diag}\, }

\begin{document}
\title{On spectral synthesis in varieties containing the solutions of inhomogeneous linear functional equations}
\author{Gergely Kiss \\ {gergely.kiss@uni.lu} \and Csaba Vincze \\ {csvincze@science.unideb.hu}}

\footnotetext[2]{{\bf MR subject classification:} primary 43A45, 43A70,
secondary 13F20}
\footnotetext[3]{G. Kiss is supported by the Internal Research Project R-STR-1041-00-Z of the University of Luxembourgh and by the Hungarian National Foundation for Scientific Research,
Grant No. K104178. Cs. Vincze is supported by the University of Debrecen's internal research project RH/885/2013.}

\maketitle

\abstract{As a continuation of our previous work \cite{KV2} the aim of the recent paper is to investigate the solutions of special inhomogeneous linear functional equations by using spectral synthesis in  translation invariant closed linear subspaces of additive/multiadditive functions containing the restrictions of the solutions to finitely generated fields. The idea is based on the fundamental work of M. Laczkovich and G. Kiss \cite{KL}. Using spectral analysis in some related varieties we can prove the existence of special solutions (automorphisms) of the functional equation but the spectral synthesis allows us to describe the entire space of solutions on a large class of finitely generated fields. It is spanned by the so-called exponential monomials which can be given in terms of automorphisms of $\cc$ and differential operators. We apply the general theory to some inhomogeneous problems motivated by quadrature rules of approximate integration \cite{KKSZ08}, see also \cite{KKSZ} and \cite{KKSZW}.

\section{Introduction and preliminaries}\label{s1}

Let $\cc$ denote the field of complex numbers. We are going to investigate the family of functional equations of type
\begin{equation}\label{ek}
\suminf=c_p\cdot \sum_{l=0}^p x^l y^{p-l} \qquad (x,y\in \cc ),
\end{equation}
where $\al_i$ and $\be_i$ $(i=1, \ldots, n)$ are given real or complex parameters, $p=1, \ldots, 2n-1$ and $c_p\in \cc$ is a constant depending on $p$. For some technical reasons we pay a special attention to the case of $p=1$, i.e.  
\begin{equation}\label{e2}
\suminf=c\cdot (x+y) \qquad (x,y\in \cc ),
\end{equation}
where $c_1$ is rewritten as $c$ for the sake of simplicity. The problem of solving the family of equations \eqref{ek} as $p$ runs through its possible values $1, \ldots, 2n-1$ is equivalent to the solution of equation
\begin{equation}\label{e1}
F(y)-F(x)=(y-x)\suminf,
\end{equation}
where $x,y\in \cc$ and $f,F:\cc\to \cc$ are unknown functions. It is motivated by quadrature rules of approximate integration \cite{KKSZ08}, see also \cite{KKSZ} and \cite{KKSZW}. 

\begin{remark}\rm{In order to substitute $x=0$ or $y=0$ into \eqref{ek} we agree that $0^0:=1$. Such a special choice reproduces the pair of equations
\begin{equation}\label{ekxy}
\sumina f(\al_i x)=c_p\cdot x^p\ \ \textrm{and}\ \ \sumina f(\be_i y)=c_p\cdot y^p \qquad (x,y\in \cc )
\end{equation}
that are consequences of \eqref{e1} for monomial solutions of degree $p$. For a more detailed survey of the preliminary results see \cite{KV2}.}
\end{remark}

Let $(G, *)$ be an Abelian group;  $\cc ^G$ denotes the set of complex valued functions defined on $G$. A function $f:G\to \cc$ is a {\it
generalized polynomial}, if there is a non-negative integer $p$ such that
\begin{equation}\label{e4}
\De _{g_1}\ldots \De _{g_{p+1}}f=0
\end{equation}
for any $g_1 \stb g_{p+1} \in G$. Here $\De _g$ is the difference operator
defined by $\De _g f(x)=f(g*x)-f(x) \ (x\in G)$, where $f\in \cc ^G$
and $g\in G$. The smallest $p$ for which \eqref{e4} holds for any $g_1 \stb g_{p+1} \in G$ is the {\it degree} of the
generalized polynomial $f$. A function $F:G^p \to \cc$ is {\it $p$-additive}, if it is additive
in each of its variables. A function $f\in \cc ^G$ is called a
{\it generalized monomial of degree $p$}, if there is a
symmetric $p$-additive function
$F$ such that $f(x)=F(x\stb x)$ for any $x\in G$. It is known that any generalized polynomial function can be written as the sum of generalized monomials \cite{SZ3}. By a general result of M. Sablik \cite{S} any solution of \eqref{e1} is a generalized polynomial of degree at most $2n-1$ under some mild conditions for the parameters in the functional equation:
\begin{enumerate}
\item $\al_1\stb \al_n, \be_1 \stb \be_n\in \rr$ or $\ \cc$,
\item $\al_i+\be_i\ne 0$,
\item \begin{equation}\label{efelt1}\left|
\begin{array}{ll}
\al_i & \be_i \\
\al_j & \be_j
\end{array} \right| \ne 0, ~~i\ne j, ~ i, j\in \{1\stb n\};\end{equation}
\end{enumerate}
see also Lemma 2 in \cite{KKSZ}. The generalized polynomial solutions of \eqref{e1} are constituted by the sum of the diagonalizations of $p$-additive functions satisfying equations of type \eqref{ek}. 

In what follows we are going to use spectral synthesis in translation invariant closed linear subspaces of additive functions on some finitely generated fields containing the restrictions of the solutions of functional equation \eqref{e2}. Note that the translation invariance is taken with respect to the multiplicative group structure. We will use spectral synthesis in some related varieties of equation \eqref{ek} for any $p>1$ too. The idea is based on \cite{KL}. To describe the space of the solutions of the inhomogeneous equation we need a non-zero particular solution in the first step. In some special cases it is enough to use spectral analysis to find such a solution \cite{KV2}. Otherwise the spectral analysis proves the existence of special solutions (automorphisms) of the functional equations in the homogeneous case; see e.g. \cite{KV}, \cite{KVV}, \cite{VV0} and \cite{KV2}. This means that we have only some necessary conditions for the existence of a nonzero solution of the inhomogeneous problem and we need the application of spectral synthesis in the varieties to give the description of the solution space on a large class of finitely generated fields. It is spanned by the so-called exponential monomials which can be given in terms of automorphisms of $\cc$ and differential operators. Unfortunately, the description of all exponential monomials spanning the entire space of the solutions seems to be beyond hope in general; see Example 2 in subsection 4.2. Our results give an explicite and unified technic to solve the problem of finding solutions at all: it is based on the spectral analysis in the first part \cite{KV2} of the investagations and the present paper completes the solution of the problem by the application of spectral synthesis.    

%






\subsection{Varieties generated by non-trivial solutions of linear functional equations}

The varieties we are going to investigate have been constructed in our previous work \cite{KV2} for the application of the so-called spectral analysis. In what follows we summarize the basic steps of the constructions. Let $G$ be an Abelian group. By a {\it variety} we mean a translation invariant closed linear
subspace of $\cc^G$.

\subsubsection{Varieties of additive solutions}

Let a finitely generated subfield $K\subset \cc$ containing the
parameters $\al_i ,\be_i$ $(i=1\stb n)$ be fixed. If $V_1$ is the set of additive functions on $K$ then it is a closed linear subspace in $\cc ^{K}$; for the proof see \cite{KL}.  

\begin{defin}
Let $S_1$ be the subset of $\ V_1$, where $\tilde{f}\in S_1$ if and only if there exists $\tilde{c}\in \cc$ such that  
\begin{equation}\label{e2tilde}
\suminftilde=\tilde{c}\cdot (x+y) \qquad (x,y\in K ).
\end{equation}
\end{defin}

By Lemma 2.3 in \cite{KV2}, $S_1$ is a closed linear subspace of $\ V_1$. Let $K^*= \{x\in K: x\neq 0\}$ be the Abelian group
with respect to the multiplication in $K$. We also put $V_1^* =\{ f|_{K^*} : f\in V_1 \}$ and $S_1^* =\{ \tilde{f}|_{K^*} : \tilde{f}\in S_1 \} .$ By Lemma 2.5 in \cite{KV2}, $V_1^{*}$ and $S_1^*$ are varieties in $\cc^{K^*}$. Recall that the translation invariance is taken with respect to the multiplicative group structure, i.e. if $\tilde{f}\in S^*_1$, then the map $\tau_a \tilde{f}\colon x\in K^* \mapsto \tilde{f}(a x)$ also belongs to $S_1^*$ for every  $a \in  K^*$. 

\begin{defin} $S_1^0$ is the subspace of $S_1$ belonging to the homogeneous case $\tilde{c}=0$.
\end{defin}

\subsubsection{Varieties generated by higher order monomial solutions}

Let a finitely generated subfield $K\subset \cc$ containing the
parameters $\al_i ,\be_i$ $(i=1\stb n)$ be fixed. If $V_p$ is the set of $p$-additive functions on $K$ then it is a closed linear subspace in $\cc^{G}$, where $G=K\times \ldots \times K$ is the Cartesian product of $K$ with itself ($p$-times); for the proof see \cite{KL}. For any $p$-additive function $F_p$ let us define $F_p^{\sigma}$ as
$$F_p^{\sigma}(w_1, \ldots, w_p):=F_p(w_{\sigma(1)}, \ldots, w_{\sigma(p)}),$$
where $\sigma$ is a permutation of the elements $1, \ldots, p$.

\begin{defin} Let $S_p$ be the subset of $\ V_p$, where $\tilde{F}_p\in S_p$ if and only if there exists $\tilde{c}_p\in \cc$ such that  
{\footnotesize{\begin{equation}\label{multimultitilde}
\begin{split}
&\sum_{i=1}^n a_i\tilde{F}_p^{\sigma}(\al_ix_1, \ldots, \al_i x_p)=\tilde{c}_p\cdot x_1\cdot \ldots \cdot x_p,\\
&\sum_{i=1}^n a_i\tilde{F}_p^{\sigma}(\be_i y_1, \ldots, \be_i y_p)=\tilde{c}_p\cdot y_1 \cdot \ldots \cdot y_p,\\
&\sum_{i=1}^n a_i \binom{p}{l} \tilde{F}_{p}^{\sigma} (\al_i x_1, \ldots, \al_i x_l, \be_i y_1, \ldots, \be_i y_{p-l})=\tilde{c}_p\cdot x_1\cdot \ldots \cdot x_l\cdot y_1 \cdot \ldots \cdot y_{p-l}
\end{split}
\end{equation}}}

\noindent
$(l=1, \ldots, p-1)$ for any permutation $\sigma$ of the elements $1, \ldots, p$.
\end{defin}

$S_p$ is a closed linear subspace of $\ V_p$. If the sets $V_p^*$ and $S_p^*$ consist of the restrictions of the elements in $V_p$ and $S_p$ to $G^*:=K^*\times \ldots \times K^*\ \ (p-\textrm{times}),$
respectively then they are varieties in $\cc^{G^*}$. Lemma 2.8 in \cite{KV2} shows that $S_p^*$ is the variety in $\cc^{G^*}$ generated by the restriction of symmetric $p$-additive functions to $G^*$ provided that the diagonalizations are the solutions of functional equation 
\begin{equation}\label{klon1}
\sum_{i=1}^n a_i f(\al_ix+\be_iy)=\tilde{c}_p\cdot \sum_{l=0}^p x^l y^{p-l} \qquad (x,y\in K ).
\end{equation}
for some $\tilde{c}_p\in \cc$.} 

\begin{defin} $S_p^0$ is the subspace of $S_p$ belonging to the homogeneous case $\tilde{c}_p=0$.
\end{defin}

\subsection{Applications of spectral analysis}

Let $(G, *)$ be an Abelian group. A function $m\colon G\to \cc$ is called \emph{exponential} if it is multiplicative: $m(x*y)=m(x)m(y)$ for any $x, y\in G$. If a variety contains an exponential function  
then we say that {\it
spectral analysis holds in this variety}. If spectral analysis holds
in each variety on $G$, then {\it spectral analysis
holds on $G$.} The main results of \cite{KV2} (Theorem 3.3 and Theorem 4.1) are based on the application of spectral analysis in the variety $S_p^*$ ($p=1, \ldots, 2n-1$). They can be summarized as follows\footnote{Since the variety $S_1^*$ contains the restrictions of additive functions, the exponential property results in an automorphism of $\cc$ by an extension process.}. 

\begin{theorem}\label{t1o1}
The existence of a nonzero additive solution of \eqref{e2} implies that there exist a finitely generated subfield $K\subset \cc$ containing $\al_i$ and $\be_i$ $(i=1, \ldots, n)$ and an automorphism $\phi\colon\mathbb{C}\to\mathbb{C}$ as the extension of an exponential element in $S_1^*$ such that 
\begin{equation}\label{eA1}
\sum_{i=1}^n a_i\phi(\al_ix+\be_iy)=\tilde{c}\cdot (x+y)\ \ (x, y\in K)
\end{equation}
for some $\tilde{c}\in \cc$. Especially,
\begin{equation}\label{eA}
\sum_{i=1}^n a_i\phi(\al_i)=\sum_{i=1}^n a_i\phi (\be_i)=\tilde{c}.
\end{equation}
If $\tilde{c}=0$ then 
\begin{equation}\label{eA0}
\sum_{i=1}^n a_i\phi(\al_ix+\be_iy)=0 \ \ (x, y\in \cc),
\end{equation}
i.e. $\phi$ is the solution of the homogeneous equation on $\cc$. If $\tilde{c}\neq 0$ then $\phi(x)=x$ $(x\in K)$ and
\begin{equation}
\label{simplified}
\sum_{i=1}^n a_i\al_i=\sum_{i=1}^n a_i\be_i=\tilde{c}\neq 0.
\end{equation}
Conversely, if \eqref{simplified} holds then $f:=(c/\tilde{c})\cdot x$ is a nonzero particular additive solution of \eqref{e2} on $\cc$.
\end{theorem}

The result says that if there are no automorphisms satisfying (\ref{eA}) with $\tilde{c}=0$, i.e. $S_1^0$ is trivial\footnote{
The so-called characteristic polynomial method helps us to investigate such an existence problem in terms of polynomials whose coefficients depend algebraically on the parameters $\al_i$ and $\be_i$ $(i=1, \ldots, n)$; \cite{V1}, see also \cite{VV-1}, \cite{V2} and \cite{VV}.} for any finitely generated subfield $K\subset \cc$ containing the parameters  $\al_i$ and $\be_i$ $(i=1, \ldots, n)$ then the only nonzero additive solution of \eqref{e2} must be the proportional of the identity function provided that \eqref{simplified} holds. In what follows we are interested in another possible case: if $\tilde{c}=0$ for any exponential function in $S_1^*$ then the exponentials give only translation parts in the solution of the inhomogeneous equation on $K$ and we need to apply spectral synthesis in the variety $S_1^*$ to decide the existence of a nonzero particular solution of the inhomogeneous equation on finitely generated fields containing the parameters $\al_i$ and $\be_i$ $(i=1, \ldots, n)$. The higher order analogue of Theorem \ref{t1o1} formulates the consequences of the application of spectral analysis in $S_p^*$: the existence of a nonzero monomial solution of degree $p>1$ of \eqref{ek} implies that there exist a finitely generated subfield $K\subset \cc$ containing $\al_i$ and $\be_i$ $(i=1, \ldots, n)$ and some automorphisms $\phi_i\colon \cc\to \cc$ $(i=1, \ldots, p)$ such that 
\begin{equation}\label{ektildeplus}
\sum_{i=1}^n a_i \diag \phi (\al_ix+\be_iy)=\tilde{c}_p\cdot \sum_{l=0}^p x^l y^{p-l} \qquad (x,y\in K)
\end{equation}
for some $\tilde{c}_p\in \cc$, where the product $\phi=\phi_1\cdot \ldots \cdot \phi_p$ is an exponential function in $S_p^*$ and $\diag$ means the diagonalization of the mappings. Especially
\begin{equation}\label{eAk}
\begin{split}
&\sum_{i=1}^n a_i\phi_1(\al_i)\cdot \ldots \cdot \phi_p(\al_i)=\sum_{i=1}^n a_i\phi_1(\be_i)\cdot \ldots \cdot \phi_p(\be_i)=\\
&\sum_{i=1}^n a_i \binom{p}{l} \phi_{\sigma(1)}(\al_i)\cdot \ldots \cdot \phi_{\sigma(l)}(\al_i)\cdot \phi_{\sigma(l+1)}(\be_i)\cdot \ldots \cdot \phi_{\sigma(p)}(\be_i)= \tilde{c}_p,
\end{split}
\end{equation}
where $l=1, \ldots, p-1$ and $\sigma$ is an arbitrary permutation of the indices. If $\tilde{c}_p=0$ then
\begin{equation}\label{ektilde1}
\sum_{i=1}^n a_i \diag \phi (\al_ix+\be_iy)=0 \qquad (x, y\in \cc),
\end{equation}
i.e. $\diag \phi$ is the solution of the homogeneous equation on $\cc$. If $\tilde{c}_p\neq 0$ then $\ph_1(x)=\ldots=\phi_p(x)=x$ $(x\in K)$ and
\begin{equation}
\label{simplified1}
\sum_{i=1}^n a_i \al_i^p=\sum_{i=1}^n a_i \be_i^p=\sum_{i=1}^n a_i\binom{p}{l} \al_i^l \be_i^{p-l}=\tilde{c}_p\neq 0\ \ (l=1, \ldots, p-1).
\end{equation}
Conversely, if \eqref{simplified1} holds then $f(x):=(c/\tilde{c}_p)\cdot x^p$ is a nonzero particular monomial solution of degree $p$ of \eqref{ek} on $\cc$; see Theorem 4.1 in \cite{KV2}. In what follows we are interested in the application of spectral synthesis in the variety $S_p^*$. It is spanned by the so-called exponential monomials which can be given in terms of automorphisms of $\cc$ and differential operators. The conditions for the exponential monomial solutions are formulated on a finitely generated field $K\subset \cc$ containing the parameters $\al_i$ and $\be_i$ $(i=1, \ldots, n)$. Therefore both the exponential function and the differential operator solutions provide solutions of the functional equation on $\cc$ by some (tipically transfinite) extension processes. Note that the direct application of the spectral theory (analysis and synthesis) to $\cc$ is impossible but any solution on $\cc$ can be naturally embedded in the varieties by a simple restriction. Therefore we also restrict our investigations to finitely generated fields to provide the area for the application of spectral theory (analysis and synthesis). 

\section{Spectral synthesis in the variety containing additive functions}

An \emph{exponential monomial} is the product of a generalized polynomial and an exponential function. If a variety $V\subset \cc^{G}$ is spanned by exponential monomials belonging to
$V$ then we say that {\it spectral synthesis holds in the variety $V$}.
If spectral synthesis holds in each variety on $G$, then {\it spectral synthesis holds on $G$.} If spectral synthesis holds in a variety $V$ then spectral analysis
holds in $V$, as well. Lemma 2.2 in \cite{SZ5} states the explicite result as follows. 

\begin{lemma}
Let $p$ be a nonzero generalized polynomial, $m$ is an exponential function on the Abelian group $G$; if the exponential monomial $p\cdot m$ belongs to the variety $V\subset \cc^G$ then $\Delta_h p \cdot m$ also belongs to $V$; especially $V$ contains the exponential function $m$ after finitely many steps of applying the difference operator to the polynomial term. 
\end{lemma}

The following result was proved in \cite{KL}. 

\begin{theorem}\label{tsy4}
Suppose that the transcendence degree of the field $K$ over $\qq$ is
finite. Then spectral synthesis holds in every variety on $K^*$
consisting of additive functions with respect to addition.
\end{theorem}

To give a more precise description of the solutions of functional equation \eqref{e2} we also need the notion of differential operators. 

\subsection{Differential operators on a finitely generated field $K$}
\label{ssapp1} Suppose that the complex numbers $t_1 \stb t_n$ are
algebraically independent over $\qq$. The elements of the field
$\qq(t_1 \stb t_n )$ are the rational functions of $t_1 \stb t_n$
with rational coefficients. By a {\it differential operator on
$\qq(t_1 \stb t_n )$} we mean an operator of the form
\begin{equation}\label{euj1}
D=\sum c_{i_1 \ldots i_n} \cd \frac{\partial^{i_1+\dots+i_n}}{\partial
t_1^{i_1}\cdots \partial t_n^{i_n}},
\end{equation}
where $\del /\del t_i$ is the usual partial derivative, the sum is
finite, the coefficient is a complex number in each term and the
exponents $i_1 \stb i_n$ are nonnegative integers. If $i_1 =\ldots
=i_n =0$, then $\partial^{i_1+\dots+i_n} /\partial t_1^{i_1}
\cdots \partial t_n^{i_n}$ means the identity operator on $\qq(t_1
\stb t_n )$. The degree of the differential operator $D$ is the
maximum of the numbers $i_1 +\ldots +i_n$ such that $c_{i_1 \ldots i_n} \ne 0$. It is clear that $\del /\del t_i$ is a derivation on $\qq (t_1
\stb t_n)$ for every $i=1\stb n$, i.e. it is an additive function satisfying the Leibniz rule. Therefore, any differential
operator on $\qq (t_1 \stb t_n)$ is the complex linear combination of finitely many maps of the form $d_1 \circ
\ldots \circ  d_k$, where $d_1 \stb d_k$ are derivations on $\qq
(t_1 \stb t_n).$ This observation motivates the following
definition.

\begin{defin} \label{duj1}
{\rm Let $K$ be a subfield of $\ \cc$. We say that the map $D:K\to
\cc$ is a} {\it differential operator} {\rm on $K$, if $D$ is the
complex linear combination of finitely many maps
of the form $d_1 \circ  \ldots \circ  d_k$, where $d_1 \stb d_k$ are
derivations on $K$. If $k=0$ then $d_1 \circ  \ldots
\circ  d_k$ means the identity function on $K$.}
\end{defin}

Differential operators in the sense of formula \eqref{euj1} and Definition \ref{duj1} mean the same objects on $K=\qq (t_1 \stb t_n)$ as the following Proposition shows; see \cite{KL}. 

\begin{proposition}\label{puj1}
Let $K$ be a subfield of $\cc$ and suppose that the elements $t_1
\stb t_n \in K$ are algebraically independent over $\qq$. If $D$ is
a differential operator on $K$ then the restriction of $D$ to $\qq (t_1 \stb t_n)$ is of the form
\eqref{euj1}.
\end{proposition}

\begin{defin}\label{action}
The action of a field automorphism $\phi\colon \cc \to \cc$ on the differential operator
\begin{equation}
D=\sum c_{i_1 \ldots i_n} \cd \frac{\partial^{i_1+\dots+i_n}}{\partial
t_1^{i_1}\cdots \partial t_n^{i_n}}
\end{equation}
is defined as 
\begin{equation}\label{actionformula}
D^{\phi}=\sum c_{i_1 \ldots i_n}' \cd \frac{\partial^{i_1+\dots+i_n}}{\partial
t_1^{i_1}\cdots \partial t_n^{i_n}},
\end{equation}
where $c_{i_1 \ldots i_n}':=\phi(c_{i_1 \ldots i_n})$ for any coefficient of $D$.
\end{defin}

\begin{remark}{\rm Note that if $L\su K\su \cc$ are fields and $D$ is a differential
operator on $L$, then $D$ can be extended to $K$ as a differential
operator. This is clear from the fact that every derivation can be
extended from $L$ to $K$. If $K$ is the algebraic
extension of $L$, then the extension of the differential operator is uniquely determined.}
\end{remark}

\section{Applications of spectral synthesis for additive solutions of linear functional equations}

Now we are going to apply spectral synthesis to the functional equation
\begin{equation}\label{e2K}
\suminf=c\cdot (x+y) \qquad (x,y\in K)
\end{equation}
and the related variety $S_1^*$. Recall that $K\subset \cc$ is a finitely generated subfield containing the
parameters $\al_i ,\be_i$ $(i=1\stb n)$, $V_1$ is the set of additive functions on $K$, $S_1\subset V_1$, where $\tilde{f}\in S_1$ if and only if there exists $\tilde{c}\in \cc$ such that  
\begin{equation}\label{e2tildeuj}
\suminftilde=\tilde{c}\cdot (x+y) \qquad (x,y\in K)
\end{equation}
and $S_1^*$ is the variety containing the restrictions of the elements in $S_1$ to the multiplicative subgroup $K^*$ of $K$. Since $K$ is a finitely generated field we can apply Theorem \ref{tsy4} to conclude that spectral synthesis holds in $S_1^*$. Therefore it is spanned by exponential monomials. By Theorem 4.2 in \cite{KL} we have the following basic theorem. 

\begin{theorem}
\label{t5}
The function $\tilde{f}\in S_1^*$ is an exponential monomial on
$K^*$ if and only if $\tilde{f}=\phi \circ D|_{K^*}$, where $\phi\colon \cc \to \cc$ is the extension of an
exponential function in $S_1^*$ to an automorphism of $\ \cc$ and $D$ is a differential operator on $K$.
\end{theorem}

In what follows we are going to test the most simple generating elements: if there is a solution of the form $\phi \circ D$ then it is also
a polynomial exponential function of the form $p\cdot m$ by Lemma 4.2 in \cite{KL}. As Lemma 2.2 in \cite{SZ5} shows, if the generalized polynomial $p$ is not constant, then $\Delta_h p\cdot m$ is
also a solution. In the most simple case $p$ is additive and the polynomial exponential function $p\cdot m$ can be written into the
form $\phi\circ d$ where $d\colon K\to K$ is a derivation; the details can be
found in \cite{KG}.

\begin{proposition}\label{lder}
Suppose that 
\begin{equation}\label{e1der}
\sumina \tilde{f}(\al_ix+\be_i y)=\tilde{c} \cdot (x+y)\ \ (x, y\in K)
\end{equation}
has a nonzero exponential monomial solution of the form $\phi\circ d$, where $\phi\colon \cc \to \cc$ is the extension of an
exponential function in $S_1^*$ to an automorphism of $\ \cc$ and
$d:K \to K$ is a derivation. If $\tilde{c}\neq 0$ then $\phi(x)=x$ $(x\in K)$,
\begin{equation}\label{e1f}
\sumina \al_i=\sumina \be_i=0
\end{equation}
and
\begin{equation}\label{e2f} \sumina d(\al_i)=\sumina d(\be_i)=\tilde{c}\neq 0.
\end{equation}
Conversely, if $\tilde{c}\neq 0$ then \eqref{e1f} and \eqref{e2f} imply that $\tilde{f}=(c/\tilde{c})\cdot d$ is a nonzero particular additive solution of \eqref{e2K}. 
\end{proposition}

\proof If $y=0$ then 
\begin{equation}
\label{e1derelso}
\sumina \tilde{f}(\al_ix)=\tilde{c} \cdot x,
\end{equation}
where $\tilde{f}$ is of the form $\phi\circ d$ and we have that
$$\sum_{i=1}^n a_i\phi(d(\alpha_ix))=\tilde{c}\cdot x,\ \ \textrm{i.e.} \ \ \sum_{i=1}^n a_i\phi(d(\alpha_i)x+\alpha_id(x))=\tilde{c}\cdot x$$
\begin{equation}
\label{zeroconclusion}
\left(\sum_{i=1}^n a_i\phi(d(\alpha_i))\right)\phi(x)+\left(\sum_{i=1}^n a_i\phi(\alpha_i)\right)\phi(dx)=\tilde{c}\cdot x.\end{equation}
If $x=1$ then
\begin{equation}
\label{firstconclusion}
\sum_{i=1}^n a_i\phi(d(\alpha_i))=\tilde{c},\ \ \textrm{i.e.}\ \ \tilde{c}\cdot \phi(x)+\left(\sum_{i=1}^n a_i\phi(\alpha_i)\right) \phi(d x)=\tilde{c}\cdot x
\end{equation}
which means that
\begin{equation}
\label{secondconclusion}
\left(\sum_{i=1}^n a_i\phi(\alpha_i)\right) \phi(d x)=\tilde{c}\cdot(x-\phi(x)).
\end{equation}
In a similar way,  
\begin{equation}\label{thirdconclusion}
\begin{split}
&\left(\sum_{i=1}^n a_i\phi(\alpha_i)\right) \phi(d y)=\tilde{c}\cdot(y-\phi(y)),\\
&\left(\sum_{i=1}^n a_i\phi(\alpha_i)\right) \phi(d (xy))=\tilde{c}\cdot(xy-\phi(xy)).
\end{split}
\end{equation}
Expanding both sides of the second equation in \eqref{thirdconclusion}, equations \eqref{firstconclusion} - \eqref{thirdconclusion} imply that $\displaystyle{\phi(x)=x}$ in case of $\tilde{c}\neq 0$.  On the other hand
$$\sum_{i=1}^n a_i \alpha_i=0 \ \ \textrm{and}\ \  \sum_{i=1}^n a_id(\alpha_i)=\tilde{c}$$
in the sense of \eqref{firstconclusion}. The corresponding relations for $\be_i$'s can be derived in a similar way by substitution $x=0$. The converse of the statement is trivial. 
\hfill $\square$

\begin{remark}\label{separation}{\rm The proof of the previous theorem shows that the result and its consequences can be formulated by separating the terms containing $x$ and $y$, respectively.}
\end{remark}

\begin{cor} If $\tilde{c}\neq 0$ for an exponential monomial $\phi\circ d$ in $S_p^*$ then the space of the additive solutions of equation \eqref{e2K} is
$$\frac{c}{\tilde{c}}\cdot d+S_1^{0},$$
where $S_1^0\subset S_1$ is the subspace belonging to the homogeneous case $\tilde{c}=0$.
\end{cor}

\begin{remark}{\rm
If $\tilde{c}\neq 0$ then equation \eqref{e1f} shows that any exponential function in $S_1^*$ solves the homogeneous equation because 
$$\sumina \phi(\alpha_i)=\sumina \phi(\be_i)\neq 0$$
gives a contradiction in the sense of Theorem \ref{t1o1}. In case of $\tilde{c}=0$ 
\begin{equation}\label{hom01}
\begin{split}
&\sumina \phi(\alpha_i)=\sumina \phi(\be_i)=0,\\
&\sumina \phi(d(\alpha_i))=\sumina \phi(d(\be_i))=0,
\end{split}
\end{equation}
i.e. both the exponential function $\phi$ and the exponential monomial $\phi\circ d$ are the solutions of the homogeneous equation. It is useful to apply the inverse automorphism on both sides to unify \eqref{hom01} for the application of the characteristic polynomial method:
\begin{equation}\label{hom02}
\begin{split}
&\sum_{i=1}^n \phi^{-1}(a_i)\al_i=\sum_{i=1}^n \phi^{-1}(a_i)\be_i=0,\\
&\sum_{i=1}^n \phi^{-1}(a_i) d(\alpha_i)=\sum_{i=1}^n \phi^{-1}(a_i)d(\be_i)=0;
\end{split}
\end{equation}
\cite{V1}, see also \cite{VV-1}, \cite{V2}, \cite{VV}.}
\end{remark}

\begin{cor}
A derivation $d\colon K\to K$ is a solution of \eqref{e1der} on $K$ with $\tilde{c}\neq 0$ if and only if 
\begin{equation}\label{e1fuj}
\sumina \al_i=\sumina \be_i=0
\end{equation}
and
\begin{equation}\label{e2fuj} \sumina d(\al_i)=\sumina d(\be_i)=\tilde{c}\neq 0.
\end{equation}
\end{cor}

\subsection{An observation on differential operators and field isomorphisms} Before testing the generating element $\phi\circ D$ of $S_1^*$ we need the following key lemma. Let $t_1 \stb
t_k$ be an algebraically independent system. For the sake of simplicity let us introduce the following abbreviations:
$$\partial_1=\frac{\partial}{\partial t_1}, \ldots, \partial_k=\frac{\partial}{\partial t_k}.$$
In the sense of Proposition \ref{puj1} any difference operator on
$\qq(t_1\stb t_k)$ is of the form
$$D:=\sum_{j_1=0}^{J_1} \ldots \sum_{j_k=0}^{J_k} c_{j_1\ldots j_k}\partial_1^{j_1} \ldots \partial_{k}^{j_k}$$
and it has a unique extension to the algebraic closure of $\qq(t_1\stb t_k)$.

\begin{lemma}\label{lprac}
Suppose that 
$$D=\sum_{j_1=0}^{J_1} \ldots \sum_{j_k=0}^{J_k} c_{j_1\ldots j_k}\partial_1^{j_1} \ldots \partial_{k}^{j_k}$$
is a differential operator on $\qq(t_1\stb t_k)$ with the uniquely determined extension to the field $K$, where $\qq(t_1\stb t_k)\subset K$ and $K$ is contained in the algebraic closure of $\qq(t_1\stb t_k)$. If
$$D(x)=c\cdot \phi(x)$$ 
for any $x \in K$, where $c$ is a nonzero complex number, $\phi \colon \cc\to \cc$ is an automorphism then 
$c_{j_1 \ldots j_k}=0$ if $j_1+\ldots+j_m\geq 1$, $c_{0 \ldots 0}=c$ and $\phi(x)=x$ $(x\in K)$.
\end{lemma}

\proof
First we show the statement for $k=1$. For the sake of simplicity let $t_1=t$ and $\partial_1=\partial$. We have that
\begin{equation}
\sum_{j=1}^{J}c_{j}\partial^{j}(x)+c_0 x= c \cdot \phi(x).
\end{equation}
If $x=1$ then $c_0=c$. Substituting $x=t^j$, $j=M\stb N$: 
\begin{eqnarray*}
\frac{M!}{(M-J)!}\cdot c_J\cdot t^{M-J}+ \dots + M\cdot c_1\cdot t^{M-1}+c\cdot t^{M}&=&c\cdot v^M,\\
\ldots&=&\ldots\\
\frac{N!}{(N-J)!}\cdot c_J\cdot t^{N-J}+ \dots + N\cdot c_1 \cdot t^{N-1}+c\cdot t^{N}&=&c\cdot v^N,
\end{eqnarray*}
where $v= \phi(t)$. Let us divide the equations by $t^M\stb t^N $, respectively:  
\begin{eqnarray*}
\frac{M!}{(M-J)!}\cdot c_J\cdot t^{-J}+ \dots + M \cdot c_1 \cdot t^{-1}+c&=&c\cdot \left(\frac{v}{t}\right)^M,\\
\ldots&=&\ldots\\
\frac{N!}{(N-L)!}\cdot c_J\cdot t^{-J}+ \dots + N\cdot c_1\cdot t^{-1}+c&=&c\cdot \left(\frac{v}{t}\right)^N
\end{eqnarray*}
Using the notation $c_j\cdot t^{-j}=\mu_j$ we can write 
\begin{eqnarray*}
\frac{M!}{(M-J)!}\cdot \mu_J+ \dots + M \cdot \mu_1 + c &=&c\cdot \left(\frac{v}{t}\right)^M,\\
\ldots&=&\ldots\\
\frac{N!}{(N-J)!}\cdot \mu_J + \dots + N\cdot \mu_1+c&=&c\cdot \left(\frac{v}{t}\right)^N.
\end{eqnarray*}
The left hand sides of the equations are polynomial expressions in $M$, $\ldots$, $N$, respectively (the degree of the polynomial is the degree $J$ of the differential operator). The right hand sides of these equations are exponential expressions in $M$, $\ldots$, $N$, repectively. It follows that $c$ must be zero or $v=t$ and $\displaystyle{v/t=1}$. In both cases we get that $$
 \begin{bmatrix} \frac{M!}{(M-J)!} & \cdots &M\\ \cdots &\cdots &\cdots \\ \frac{N!}{(N-J)!}\ &\cdots &N\end{bmatrix} \times \left[ \begin{array}{c} \mu_J \\\cdots \\\mu_1 \end{array} \right]= \left[ \begin{array}{c} 0 \\ \ldots \\ 0 \end{array} \right].
$$
If $N-M=J$ then it is a quadratic matrix. Using that
$$M(M-1)+M=M^2,\ \ M(M-1)(M-2)+3M(M-1)+M=M^3, \ \ \ldots$$
we have that its determinant is the same as that of the usual Vandermonde-matrix
$$
 \begin{bmatrix} M^J & \cdots &M\\ \cdots &\cdots &\cdots\\ N^J\ &\cdots &N\end{bmatrix}.
$$
The nonzero determinant implies that $\mu_j=0$ for any $j=1 \stb J$, i.e. $c_1=\ldots=c_J=0$ and $c_0=c$ as we have seen above. Using that the extension of $D$ to $K$ is uniquely detemined it follows that $D(x)=c\cdot x$ $(x\in K)$ and, consequently, $\phi(x)=x$ $(x\in K)$. We sketch the inductive step  for $k=2$. Let $D$ be a differential operator of the variables $t_1$ and $t_2$. If one of the variables, say $t_2$, is keeping constant then we can repeat the previous procedure by the substitution of  $t_1^M\cdot t_2^{J_2}, \ldots, t_1^N\cdot t_2^{J_2}$, where $N-M=J_1$. We have that for any $j_1=1, \ldots, J_1$ 
$$\mu_{j_1}:=\sum_{j_2=0}^{J_2}c_{j_1j_2}\frac{J_2!}{(J_2-j_2)!} t_2^{J_2-j_2}\cdot t_1^{-j_1}=0,$$
i.e.
$$\sum_{j_2=0}^{J_2}c_{j_1j_2}\frac{J_2!}{(J_2-j_2)!} t_2^{J_2-j_2}=0,$$
where the left hand side is a polynomial expression of the variable $t_2$. This means that $c_{j_1 j_2}=0$, where  $j_1=1, \ldots, J_1$ and $j_2=0, \ldots, J_2$. Changing the role of the variables: $c_{j_1 j_2}=0$, where  $j_1=0, \ldots, J_1$ and $j_2=1, \ldots, J_2$. Therefore $D(x)=c\cdot x$ for any $x\in \qq(t_1, t_2)$. Using that the extension of $D$ to $K$ is uniquely detemined it follows that $D(x)=c\cdot x$ $(x\in K)$ and, consequently, $\phi(x)=x$ $(x\in K)$.
 \hfill $\square$
\bsk

\begin{remark} \rm{Lemma \ref{lprac} is of the same type as the statement in \cite{GP} about the linear independency of the iterates of any nonzero real derivation.}
\end{remark}

In what follows we are going to test the generating elements of the form $\phi\circ D$ by the transcendence degree of the parameters $\alpha_1$, $\ldots$, $\alpha_n$, $\be_1$, $\ldots$, $\be_n$ over the rationals.

\subsection{The case of transcendence degree $0$ (algebraic parameters)}

Suppose that the parameters $\alpha_1$, $\ldots$, $\alpha_n$, $\be_1$, $\ldots$, $\be_n$ are algebraic numbers over the rationals. Since there is no nontrivial derivation on $K=\qq(\al_1, \ldots, \alpha_n, \beta_1, \ldots, \beta_n)$, any element $\phi\circ D\in S_1$ reduces to $\phi$; see e.g. \cite{KG} and \cite{K}. By subsection 1.2 this means that 
$$\sumina \phi(\al_i)=\sumina \phi(\be_i)=\tilde{c}$$
and we have two possibilities: if $\tilde{c}\neq 0$ then $\phi(x)=x$ $(x\in K)$ and 
any solution of \eqref{e2K} must be of the form
\begin{equation}
\label{boing}
\frac{c}{\tilde{c}} \cdot x+ \sum_j c_j\cdot \phi_j(x) \ \ (x\in K)
\end{equation}
provided that $\sumina \al_i=\sumina \be_i=\tilde{c}\neq 0$ (Theorem \ref{t1o1}). The second term in \eqref{boing} contains the linear combination of automorphisms satisfying 
\begin{equation}\label{e2alg}
\sumina \phi_j(\al_i)=\sumina \phi_j(\be_i)=0.
\end{equation}
According to the characteristic polynomial method \cite{V1}, see also \cite{VV-1}, \cite{V2}, \cite{VV} there is only finitely many different $\phi_j$'s, i.e. the space $S_1^0$ is finitely generated. If 
\begin{equation}
\label{zerozero}
\sumina \al_i=\sumina \be_i=0
\end{equation}
then both the exponential functions and the differential operators solve the homogeneous equation\footnote{This means that they belong to $S_1^0$.}. Therefore we have no solutions of \eqref{e2K} with $c\neq 0$. 

\begin{remark} {\rm Suppose that 
$K$ is a finitely generated field containing the algebraic parameters $\al_1$, $\ldots$, $\al_n$ and $\be_1$, $\ldots$, $\be_n$ such that $K$ contains at least one transcendental number to provide the existence of nontrivial differential operators. Then we can formulate the following result too: any solution of \eqref{e2K} must be of the form   
$$\frac{c}{\tilde{c}}\cdot x+ \sum c_j\phi_j \circ D_j(x)$$
provided that $\sumina \al_i=\sumina \be_i=\tilde{c}\neq 0$, the automorphism $\phi_j$ satisfies \eqref{e2alg} and $D_j$ is an arbitrary differential operator on $K$ for any $j=1, \ldots, n$. Observe that if the automorphism part $\phi$ solves the homogeneous equation then $\phi\circ D$ is also a solution of the homogeneous equation for any differential operator $D$. The case of algebraic parameters is closely related to the theory of spectral analysis because of the trivial action of the differential operators on algebraic elements. The transcendence degree of the embedding (finitely generated) field has no any influence in this sense. In case of \eqref{zerozero} we have no solutions of \eqref{e2K} with $c\neq 0$.}
\end{remark}

\subsection{The case of higher transcendence degree}

Suppose that the transcendence degree of the parameters $\alpha_1$, $\ldots$, $\alpha_n$, $\be_1$, $\ldots$, $\be_n$ is $k$, i.e. we have a field extension $L=\mathbb{Q}(t_1, \ldots, t_k)$ such that $t_1$, $\ldots$, $t_k$ are algebraically independent over the rationals and the algebraic closure of $L$ contains the parameters $\alpha_i$'s and $\beta_i$'s. Let us introduce the abbreviations
$$\partial_1=\frac{\partial}{\partial t_1}, \ldots, \partial_k=\frac{\partial}{\partial t_k}.$$
By a simple induction we have that 
\begin{equation}
\label{eval1}
\partial^n (xy)=\sum_{k=0}^n {n \choose k} \partial^k (x) \partial^{n-k} (y)\ \ (x, y \in K),
\end{equation}
where $K$ is a finitely generated field containing the parameters $\al_i$ and $\be_i$ $(i=1, \ldots, n)$ such that $L\subset K$ and $K$ is contained in the algebraic closure of $L$. Using formula (\ref{eval1})
$$
\partial_1^{j_1} \ldots \partial_k^{j_k}(xy)=\partial_1^{j_1} \ldots \partial_{k-1}^{j_{k-1}}\left(\sum_{l_k=0}^{j_k} {j_k \choose l_k} \partial_k^{j_k-l_k}(x)\partial_k^{l_k}(y)\right)=$$
$$
\sum_{l_k=0}^{j_k} \ldots \sum_{l_1=0}^{j_1}  {j_k \choose l_k}\cdot \ldots \cdot {j_1 \choose l_1}\partial_1^{j_1-l_1} \ldots \partial_{k}^{j_{k}-l_k}(x)\partial_1^{l_1} \ldots \partial_{k}^{l_k}(y)\ \ (x, y\in K).
$$
Recall that the action of an automorphism $\phi$ on a differential operator 
$$D:=\sum_{j_1=0}^{J_1} \ldots \sum_{j_k=0}^{J_k} c_{j_1\ldots j_k}\partial_1^{j_1} \ldots \partial_{k}^{j_k}
$$
is defined as 
$$
D^{\phi}=\sum_{j_1=0}^{J_1} \ldots \sum_{j_k=0}^{J_k} c_{j_1\ldots j_k}'\partial_1^{j_1} \ldots \partial_{k}^{j_k},
$$
where $c_{i_1 \ldots i_n}':=\phi(c_{i_1 \ldots i_n})$ for any coefficient of $D$; see formula \eqref{actionformula}. In what follows we frequently need the following family of differential operators generated by $D$: if $0\leq m_1 \leq J_1$, $\ldots$, $0\leq m_k \leq J_k$ then
$$D_{m_1 \ldots m_k}:=\sum_{j_1=m_1}^{J_1} \ldots \sum_{j_k=m_k}^{J_k} c_{j_1\ldots j_k} {j_k \choose m_k}\cdot \ldots \cdot  {j_1 \choose m_1}\partial_1^{j_1-m_1} \ldots \partial_{k}^{j_{k}-m_k};$$
especially $D_{0 \ldots 0}=D.$ Therefore
$$D_{m_1 \ldots m_k}^{\phi}:=\sum_{j_1=m_1}^{J_1} \ldots \sum_{j_k=m_k}^{J_k} c_{j_1\ldots j_k}' {j_k \choose m_k}\cdot \ldots \cdot  {j_1 \choose m_1}\partial_1^{j_1-m_1} \ldots \partial_{k}^{j_{k}-m_k}.$$

\begin{proposition}\label{ldiffuj}
Suppose that 
\begin{equation}\label{e1diffuj}
\sumina \tilde{f}(\al_ix+\be_i y)=\tilde{c} \cdot (x+y)\ \ (x, y\in K)
\end{equation}
has a nonzero exponential monomial solution of the form $\phi\circ D$, where $\phi\colon \cc \to \cc$ is the extension of an
exponential function in $S_1^*$ to an automorphism of $\ \cc$ and
\begin{equation}
\label{difop1uj}
D:=\sum_{j_1=0}^{J_1} \ldots \sum_{j_k=0}^{J_k} c_{j_1\ldots j_k}\partial_1^{j_1} \ldots \partial_{k}^{j_k}
\end{equation}
is a differential operator on $K$ by its uniquely determined extension to the algebraic closure of $L$. 
If $\tilde{c}\neq 0$ then $\phi(x)=x$ $(x\in K)$, 
\begin{equation}
\label{system1diffuja}
\sum_{i=1}^n a_i D_{m_1 \ldots m_k}^{\phi}(\alpha_i)=\sum_{i=1}^n a_iD_{m_1 \ldots m_k}^{\phi}(\be_i)=0
\end{equation}
if $m_1+\ldots+m_k\geq 1$ and 
\begin{equation}
\label{system1diffujb}
\sum_{i=1}^n a_i D^{\phi}(\alpha_i)=\sum_{i=1}^n a_iD^{\phi}(\be_i)=\tilde{c}.
\end{equation}
Conversely, if $\tilde{c}\neq 0$ then \eqref{system1diffuja} and \eqref{system1diffujb} imply that  $f:=(c/\tilde{c})\cdot D^{\phi}$ is a nonzero particular additive solution of \eqref{e2K}.
\end{proposition}

\proof
If $y=0$ then equation (\ref{e1diffuj}) reduces to
\begin{equation}
\label{e1derauj}
\sumina \tilde{f}(\al_ix)=\tilde{c} \cdot x,
\end{equation}
where, by our assumption, $\tilde{f}$ is of the form $\phi\circ D$,
$$D:=\sum_{j_1=0}^{J_1} \ldots \sum_{j_k=0}^{J_k} c_{j_1\ldots j_k}\partial_1^{j_1} \ldots \partial_{k}^{j_k}.$$
Substituting in (\ref{e1derauj})
$$\sum_{i=1}^n a_i\phi(D(\alpha_i x))=\tilde{c}\cdot x,$$
$$\sum_{i=1}^n a_i \sum_{j_1=0}^{J_1} \ldots \sum_{j_k=0}^{J_k} \phi(c_{j_1\ldots j_k})\phi(\partial_1^{j_1} \ldots \partial_{k}^{j_k}(\alpha_ix))=\tilde{c}\cdot x$$
$$\sum_{i=1}^n a_i \sum_{j_1=0}^{J_1} \ldots \sum_{j_k=0}^{J_k} \phi(c_{j_1\ldots j_k})\sum_{l_k=0}^{j_k} \ldots \sum_{l_1=0}^{j_1}  {j_k \choose l_k}\cdot \ldots \cdot   {j_1 \choose l_1}$$
$$\phi(\partial_1^{j_1-l_1} \ldots \partial_{k}^{j_{k}-l_k}(\alpha_i))\phi(\partial_1^{l_1} \ldots \partial_{k}^{l_k}(x))=\tilde{c}\cdot x.$$
Applying the inverse automorphism $\phi^{-1}$ to both sides 
$$\sum_{i=1}^n \phi^{-1}(a_i) \sum_{j_1=0}^{J_1} \ldots \sum_{j_k=0}^{J_k} c_{j_1\ldots j_k}\sum_{l_k=0}^{j_k} \ldots \sum_{l_1=0}^{j_1}  {j_k \choose l_k}\cdot \ldots \cdot  {j_1 \choose l_1}$$
$$\partial_1^{j_1-l_1} \ldots \partial_{k}^{j_{k}-l_k}(\alpha_i)\partial_1^{l_1} \ldots \partial_{k}^{l_k}(x)=\phi^{-1}(\tilde{c})\phi^{-1}(x).$$
If $\tilde{c}\neq 0$ then, by Lemma \ref{lprac},
$$\lambda_{m_1 \ldots m_k}:=\sum_{i=1}^n \phi^{-1}(a_i) \sum_{j_1=m_1}^{J_1} \ldots \sum_{j_k=m_k}^{J_k} c_{j_1\ldots j_k}  {j_k \choose m_k}\cdot  \ldots \cdot {j_1 \choose m_1}$$ 
$$\partial_1^{j_1-m_1} \ldots \partial_{k}^{j_{k}-m_k}(\alpha_i)=0$$
if $m_1+\ldots+m_k\geq 1$,  
$$\lambda_{0\ldots 0}:=\sum_{i=1}^n \phi^{-1}(a_i) \sum_{j_1=0}^{J_1} \ldots \sum_{j_k=0}^{J_k} c_{j_1\ldots j_k}\partial_1^{j_1} \ldots \partial_{k}^{j_{k}}(\alpha_i)=\phi^{-1}(\tilde{c})$$
if $m_1=\ldots=m_k=0$ and $\phi^{-1}(x)=x$ ($x\in K$). Taking the action of $\phi$ on both sides of the equations it follows that  
$$\sum_{i=1}^n a_i D_{m_1 \ldots m_k}^{\phi}(\alpha_i)=0\ \ (m_1+\ldots+m_k\geq 1)\ \ \textrm{and} \ \ \sum_{i=1}^n a_i D^{\phi}(\alpha_i)=\tilde{c}$$
because of $D_{0 \ldots 0}=D$ and $\phi(x)=x$ for any $x\in K$. Note that the terms of the form  $\partial_1^{j_1-m_1} \ldots \partial_{k}^{j_{k}-m_k}(\alpha_i)$ also belongs to $K$ (see Remark \ref{computediff}).
The corresponding system of equations for the para\-meters $\be_i$'s can be derived in a similar way by substitution $x=0$. The converse of the statement is trivial. 
\hfill $\square$

\begin{cor} If $\tilde{c}\neq 0$ for an exponential monomial $\phi\circ D$ in $S_p^*$ then the space of the additive solutions of equation \eqref{e2K} is
$$\frac{c}{\tilde{c}}\cdot D^{\phi}+S_1^{0},$$
where $S_1^0\subset S_1$ is the subspace belonging to the homogeneous case $\tilde{c}=0$.
\end{cor}

\begin{remark}\label{matrixforms}{\rm
Equations \eqref{system1diffuja} and \eqref{system1diffujb} can be considered as linear systems of equations for the unknown quantities $c_{0 \ldots 0}'$, $\ldots$, $c_{j_1 \ldots j_k}'$, $\ldots$, $c_{J_1 \ldots J_k}'$ in the expression of the solution $D^{\Phi}$. Following the lexicographic ordering we have upper triangle matrices from the definition of $D_{m_1 \ldots m_k}$. The coefficients of $c_{m_1 \ldots m_k}'$ in equations \eqref{system1diffuja} and \eqref{system1diffujb} are $\displaystyle{\sum_{i=1}^n a_i \al_i}$ and $\displaystyle{\sum_{i=1}^n a_i \be_i}$, respectively. Suppose (for example) that the matrix containing $\al_i$'s is regular, i.e. 
$\displaystyle{\sum_{i=1}^n a_i \al_i\neq 0}$. For an upper triangle fundamental matrix Cramer's rule says that
$$c_{0\ldots 0}'=\tilde{c}\cdot \frac{\ \ \ \left(\sum_{i=1}^n a_i \al_i\right)^{N-1}}{\left(\sum_{i=1}^n a_i \al_i\right)^N}=\frac{\tilde{c}}{\sum_{i=1}^n a_i \al_i},$$
where $N=J_1\cdot \ldots \cdot J_k$ and $c_{j_1 \ldots j_k}=0$ if $j_1+\ldots+j_k >0$. Therefore we have two possible cases:
$$\sum_{i=1}^n a_i \al_i=\sum_{i=1}^n a_i \be_i\neq 0$$
and $D^{\phi}$ reduces to the proportional of the identity:
$$D^{\phi}(x)=c'\cdot x\ \ (x\in K),\ \ \textrm{where}\ \ c'=\frac{\tilde{c}}{\sum_{i=1}^n a_i \al_i}=\frac{\tilde{c}}{\sum_{i=1}^n a_i \be_i}.$$
Otherwise
$$\sum_{i=1}^n a_i \al_i=\sum_{i=1}^n a_i \be_i=0$$
and we can reduce the order of the upper triangle matrices of systems \eqref{system1diffuja} and \eqref{system1diffujb}: $c_{0\ldots 0}'$ can be arbitrarily choosen and we can delete the first column and the last row to give a reduced system of linear equations for the unknown quantities $c_{0 \ldots 0 1}'$, $\ldots$, $c_{j_1 \ldots j_k}'$, $\ldots$, $c_{J_1 \ldots J_k}'$. This means that if $D$ is an at least first order differential operator then 
$$\sum_{i=1}^na_i \alpha_i=\sum_{i=1}^na_i \beta_i=0$$
implies that any exponential function $\phi$ must be in $S_1^0$, i.e. $\phi$ is the solution of the homogeneous equation on $K$. By contraposition, if there is an exponential function in $S_1^*$ with $\tilde{c}\neq 0$ then it must be the proportional of the identity function on $K$ together with the exponential monomial $\phi \circ D$, i.e. $D$ is a differential operator of degree zero.}
\end{remark}

According to the technical difficulties of the discussion of the matrix forms of \eqref{system1diffuja} and \eqref{system1diffujb} we omit the details in general. The case of transcendence degree $1$ can be entirely solved by the technic of linear systems of equations as we shall see in subsection 3.3.1; see also subsection 5.2.1. In case of the higher transcendence degree we omit the further theoretical computations but we present how the method is working in explicite cases; section 4.

\begin{remark} \label{homogenity} {\rm If $\tilde{c}=0$ then 
\begin{equation}\label{diffophom}
\sum_{i=1}^n \phi^{-1}(a_i) D_{m_1 \ldots m_k}(\al_i)=\sum_{i=1}^n \phi^{-1}(a_i) D_{m_1 \ldots m_k}(\be_i)=0
\end{equation}
for any $0\leq m_i \leq J_i$ $(i=1, \ldots, k)$. These equations can be also considered as linear systems of equations for the unknown quantities $c_{0 \ldots 0}$, $\ldots$, $c_{j_1 \ldots j_k}$, $\ldots$, $c_{J_1 \ldots J_k}$; cf. Remark \ref{matrixforms}. The diagonal elements of the fundamental matrices are
$\displaystyle{\sum_{i=1}^n \phi^{-1}(a_i) \al_i}$ and $\displaystyle{\sum_{i=1}^n \phi^{-1}(a_i) \be_i}$, respectively. They have zero determinant because of
$\displaystyle{\sumina \phi(\al_i)=\sumina \phi(\be_i)=0},$ where $\phi$ is the extension of an exponential function in ${S_1^{0}}^*$ to an automorphism of $\cc$; 
cf. equation \eqref{hom02}.} 
\end{remark}

Unfortunately the action of $\phi$ is hard to compute in general. We need to apply the characteristic polynomial method to \eqref{diffophom}. This results in a system of polynomial equations. The transcendence degree $l$ of $a_i$'s over the rational gives the number of the variables of the polynomials. Their coefficients depend on 
\begin{itemize}
\item the particular actions 
$$\partial_1^{j_1-m_1} \ldots \partial_{k}^{j_{k}-m_k}(\alpha_i)\ \ \textrm{or}\ \ \partial_1^{j_1-m_1} \ldots \partial_{k}^{j_{k}-m_k}(\beta_i)$$
of the differential operator (see Remark \ref{computediff}),
\item the coefficients of the defining polynomial of an algebraic elements $u$ over $\qq(a_1, \ldots, a_l)$, where $a_1, \ldots, a_l$ form a maximal algebraically independent system and $u$ is choosen such that the missing outer parameters  belong to the simple algebraic extension $\qq(a_1, \ldots, a_l)(u)$,
\item the (rational) coefficients of the polynomials $p_j$ and $q_j$, where
$$a_j=p_j(a_1, \ldots, a_k)/q_j(a_1, \ldots, a_k)\ \ j=l+1, \ldots, n.$$ 
\end{itemize}
For the details we can refer to \cite{V1}, see also \cite{VV-1}, \cite{V2}, \cite{VV}. Practically we simultaneously determine the actions $w_1:=\phi^{-1}(a_1)$, $\ldots$, $w_l:=\phi^{-1}(a_l)$ and the quantities $c_{0 \ldots 0}$, $\ldots$, $c_{j_1 \ldots j_k}$, $\ldots$, $c_{J_1 \ldots J_k}$ as the solutions of a system of (multivariate) polynomial equations. They constitute an exponential monomial in $S_1^0$ by the composition $\phi\circ D$. 

\begin{remark}\label{separation1}{\rm The proof of the previous theorem shows that the result and its consequences can be formulated by separating the terms containing $x$ and $y$, respectively.}
\end{remark}

\begin{cor}
\label{simpleextension}
A differential operator
$$D=\sum_{j_1=0}^{J_1} \ldots \sum_{j_k=0}^{J_k} c_{j_1\ldots j_k}\partial_1^{j_1} \ldots \partial_{k}^{j_k}$$
is a solution of \eqref{e1diffuj} on $K$ with $\tilde{c}\neq 0$ if and only if 
\begin{equation}
\label{system3diffuja}
\sum_{i=1}^n a_iD_{m_1 \ldots m_k}(\alpha_i)=\sum_{i=1}^n a_i D_{m_1 \ldots m_k}(\be_i)=0
\end{equation}
if $m_1+\ldots+m_k\geq 1$ and 
\begin{equation}
\label{system3diffujb}
\sum_{i=1}^n a_i D(\alpha_i)=\sum_{i=1}^n a_i D(\be_i)=\tilde{c}\neq 0.
\end{equation}
\end{cor}

\begin{remark}\label{computediff}{\rm To compute the particular action of a differential operator on an algebraic element $u$ over $L=\mathbb{Q}(t_1, \ldots, t_k)$ we need its defining polynomial 
$$u^{m}+r_{m-1}u^{m-1}+\ldots r_1u+r_0=0.$$
Then 
$$(m u^{m-1}+(m-1)r_{m-1}u^{m-2}+...+r_1)\partial_i(u)+\sum_{j=0}^{m-1}(\partial_i r_j)u^j=0,$$
where $r_j=p_j(t_1, \ldots, t_k)/q_j(t_1, \ldots,t_k)$ is a rational fraction and $\partial_i r_j$ means the usual partial differentiation. Observe that $\partial_i(u)$ is also an algebraic number over $L=\mathbb{Q}(t_1, \ldots, t_k)$ because the algebraic numbers form a field and  $m=\textrm{deg}\ u$ implies that
$$m u^{m-1}+(m-1)r_{m-1}u^{m-2}+...+r_1\neq 0.$$
Moreover $\partial_i(u)\in \mathbb{Q}(t_1, \ldots, t_k)(u)$ and the process can be repeated to compute the action $\partial_{1}^{j_1}\ldots \partial_k^{j_k}(u)$ of the higher order term of the differential operator.}
\end{remark}

\subsubsection{The case of transcendence degree $1$}
\label{transcdegree1}

This special case gives the best result from the viewpoint of the application of the general theory. It is due to the relatively simple matrix form of the equations and the uniquely determined main terms in the differential operators. Suppose that the transcendence degree of the parameters $\alpha_1$, $\ldots$, $\alpha_n$, $\be_1$, $\ldots$, $\be_n$ is $1$, i.e. we have a field extension $L=\qq(t)$ for some transcendental number over the rationals such that its algebraic closure contains the parameters $\alpha_i$'s and $\beta_i$'s. In case of $k=1$ systems \eqref{system1diffuja} and \eqref{system1diffujb} can be written in a more detailed form (cf. Remark \ref{matrixforms}). For example 
\begin{equation}
\label{diagonal1}
\begin{split}
& c_0'\sum_{i=1}^na_i\alpha_i+c_1'\sum_{i=1}^na_i\partial(\alpha_i)+c_2'\sum_{i=1}^na_i\partial^2(\alpha_i)+\ldots+c_J'\sum_{i=1}^na_i\partial^J(\alpha_i)=\tilde{c}\\
& c_1'\sum_{i=1}^na_i\alpha_i+{2 \choose 1}c_2'\sum_{i=1}^n a_i\partial(\alpha_i)+\ldots+c_J'{J \choose 1}\sum_{i=1}^n a_i\partial^{J-1}(\alpha_i)=0\\
& ...\\
& c_{J-1}'\sum_{i=1}^na_i\alpha_i+{J \choose J-1}c_J'\sum_{i=1}^na_i\partial(\alpha_i)=0\\
& c_J'\sum_{i=1}^n a_i\alpha_i=0.
\end{split}
\end{equation}

The corresponding system of equations containing the parameters $\be_i$'s is of the same form. Since we have upper triangle matrices the following result can be easily concluded\footnote{Recall that $K$ is a finitely generated field contain\-ing the parameters $\al_i$'s and $\be_i$'s such that $L\subset K$ and $K$ is contained in the algebraic closure of $L$.}. 

\begin{proposition}\label{ldiffcor}
Suppose that 
\begin{equation}\label{e1diffcor}
\sumina \tilde{f}(\al_ix+\be_i y)=\tilde{c} \cdot (x+y)\ \ (x, y\in K) 
\end{equation}
has a nonzero exponential monomial solution of the form $\phi \circ D$ on $K$, where $\phi\colon \cc \to \cc$ is the extension of an
exponential function in $S_1^*$ to an automorphism of $\ \cc$ and
\begin{equation}
\label{difop1cor}
D=\sum_{j=0}^J c_j \partial^j, \ \ c_J\neq 0
\end{equation}
is a differential operator on $K$ by its uniquely determined extension to the algebraic closure of $L$. If $\tilde{c}\neq 0$ then $\phi(x)=x$ $(x\in K)$, 
{\footnotesize{\begin{equation}
\label{system2diffa}
\begin{split}
&\sum_{i=1}^na_i \alpha_i=0,\ \sum_{i=1}^na_i \partial(\alpha_i)=0,\ \sum_{i=1}^na_i \partial^2(\alpha_i)= 0, \ \ldots,\ \sum_{i=1}^na_i \partial^{J-1}(\alpha_i)=0,\\
&\sum_{i=1}^na_i \be_i=0,\ \sum_{i=1}^na_i \partial(\be_i)=0,\ \sum_{i=1}^na_i \partial^2(\be_i)= 0, \ \ldots,\ \sum_{i=1}^na_i \partial^{J-1}(\be_i)=0,\\
\end{split}
\end{equation}}}

\noindent
i.e. the coefficients $c_0, \ldots, c_{J-1}$ can be arbitrarily choosen and
{\footnotesize{\begin{equation}
\label{system2diffb}
\begin{split}
&\sum_{i=1}^na_i \partial^J(\alpha_i)=\sum_{i=1}^na_i \partial^J(\be_i) \neq 0,\ \ c_J'=\frac{\tilde{c}}{\sum_{i=1}^na_i \partial^J(\alpha_i)}=\frac{\tilde{c}}{\sum_{i=1}^na_i \partial^J(\be_i)},
\end{split}
\end{equation}}}

\noindent
where $c_J'=\phi(c_J).$ Conversely, if $\tilde{c}\neq 0$ then \eqref{system2diffa} and \eqref{system2diffb} imply that $f:=(c/\tilde{c})\cdot D^{\phi}$ is a nonzero particular additive solution of \eqref{e2K}, where $$D^{\phi}=\sum_{j=0}^J c_j' \partial^j, \ \ c_J'\neq 0,$$
the coefficients  $c_0', \ldots, c_{J-1}'$ can be arbitrarily choosen and $c_J'$ is determined by \eqref{system2diffb}. 
\end{proposition}

\begin{remark} {\rm If $\tilde{c}=0$ then 
\begin{equation}\label{spec1}
\sum_{i=1}^n \phi^{-1}(a_i) D_{m}(\al_i)=\sum_{i=1}^n \phi^{-1}(a_i) D_{m}(\be_i)=0
\end{equation}
for any $0\leq m \leq J$. The analogue result of \eqref{hom02} up to order $J$ of the differential operator can be easily concluded by using \eqref{diagonal1}.}
\end{remark}

\begin{cor}
\label{ldiffcoruj}
A differential operator
$$D=\sum_{j=0}^J c_j \partial^j,\ \ c_J\neq 0$$
is a solution of \eqref{e1diffcor} with $\tilde{c}\neq 0$ if and only if \eqref{system2diffa} and \eqref{system2diffb} are satisfied with $c_J=c_J'$.
\end{cor}

\begin{remark} \rm{Note that the existence of a uniquely determined main term of the differential operator is an essential difference relative to the case of higher transcendence degree; see e.g. subsection 4.2.} 
\end{remark}

\section{Examples}

The following examples illustrate how the results are working in explicite cases. For the sake of simplicity we use their separated version in the sense of Remark \ref{separation} and Remark \ref{separation1}. In case of transcendence degree $1$ we can describe the entire space of solutions (Example 1). Otherwise the situation becomes more difficult because of the missing main term of the differential operators (Example 2).

\subsection{Example 1} Let $c\neq 0$ be a complex number and consider functional equation
\begin{equation}
\label{exm02}
t^3f(tx)-t^2f(t^2x)-tf(t^3x)+f(t^4x)=c\cdot x,
\end{equation}
i.e. $n=4$ and $a_1=t^3$, $a_2=-t^2$, $a_3=-t$, $a_4=1$, $\al_1=t$, $\al_2=t^2$, $\al_3=t^3$, $\al_4=t^4$,  
where $t$ is a transcendental number over the rationals, $L=\qq(t)$ and $L\subset K$ such that $K$ is a finitely generated field and it is contained in the algebraic closure of $L$. Using that
\begin{equation}
\label{exm020}
\begin{split}
&\sum_{i=1}^4 a_i\al_i=t^3\cdot t-t^2\cdot t^2-t\cdot t^3+t^4=0\\
&\sum_{i=1}^4 a_i \partial (\al_i)=t^3-t^2\cdot(2t)-t\cdot(3t^2)+4t^3=0
\end{split}
\end{equation}
but 
\begin{equation}
\label{exm02no0}
\begin{split}
&\sum_{i=1}^4 a_i\partial^2(\al_i)=-2t^2-6t^2+12t^2=4t^2\neq 0
\end{split}
\end{equation}
we have, by Corollary \ref{ldiffcoruj}, that the differential operator
$$D=c_0+c_1\partial+\frac{c}{\ 4t^2}\partial^2\ \ (c_0, c_1 \in \cc)$$
is a non-zero particular solution of equation (\ref{exm02}) and the space of the solutions on $K$ is $D+S_1^0.$ Equation $\displaystyle{\sum_{i=1}^4 a_i\al_i}=0$ implies that any exponential function in $S_1$ belongs to $S_1^0$, i.e. 
$$t^3\phi(t)-t^2\phi(t^2)-t\phi(t^3)+\phi(t^4)=0.$$
Substituting $s=\phi(t)$ we have
$$t^3\cdot s-t^2\cdot s^2-t\cdot s^3+s^4=0\ \ \Rightarrow\ \  s\cdot(s-t)\cdot(s^2-t^2)=0.$$
Therefore $\phi(t)=t$ or $\phi(t)=-t$. To find the generating elements of the form $\phi\circ D$ in $S_1^0$ we need to use
system \eqref{spec1}:
\begin{itemize}
\item If $\phi(t)=t$ then $\phi^{-1}(a_i)=a_i$ $(i=1, \ldots, 4)$ and we have 
{\footnotesize{\begin{equation}
\begin{split}
& c_0\sum_{i=1}^4 a_i\alpha_i+c_1\sum_{i=1}^4 a_i\partial(\alpha_i)+c_2\sum_{i=1}^4 a_i\partial^2(\alpha_i)+\ldots+c_J\sum_{i=1}^4 a_i\partial^J(\alpha_i)=0\\
& c_1\sum_{i=1}^4 a_i\alpha_i+{2 \choose 1}c_2\sum_{i=1}^4 a_i\partial(\alpha_i)+\ldots+c_J{J \choose 1}\sum_{i=1}^4 a_i\partial^{J-1}(\alpha_i)=0\\
& ...\\
& c_{J-1}\sum_{i=1}^4 a_i\alpha_i+{J \choose J-1}c_J \sum_{i=1}^4 a_i\partial(\alpha_i)=0\\
& c_J \sum_{i=1}^4 a_i\alpha_i=0.
\end{split}
\end{equation}}}
\end{itemize}

\noindent
By equations \eqref{exm020} and \eqref{exm02no0} it follows that $c_2=\ldots =c_J=0$, i.e. $D=c_0+c_1\partial$ $(c_0, c_1\in \cc)$.
\begin{itemize}
\item If $\phi(t)=-t$ then $\phi^{-1}(a_1)=-a_1$, $\phi^{-1}(a_2)=a_2$, $\phi^{-1}(a_3)=-a_3$ and $\phi^{-1}(a_4)=a_4$.
System \eqref{spec1} has vanishing diagonal elements because of
$$\sum_{i=1}^4\phi^{-1}(a_i)\al_i=-t^3\cdot t-t^2\cdot t^2+t\cdot t^3+t^4=0$$
but
$$\sum_{i=1}^4\phi^{-1}(a_i) \partial(\al_i)=-t^3-t^2\cdot (2t)+t\cdot (3t^2)+4t^3=4t^3\neq 0.$$
\end{itemize}
Therefore $c_1=\ldots =c_J=0$ and $D$ reduces to the proportional of the identity function on $K$. 

This means that $S_1^0$ is spanned by the extensions of $\phi_1$, $\phi_2$ and $\phi_1\circ D$, where $\phi_1(t)=t$, $\phi_2(t)=-t$ and $D(x)=c_0 x+c_1\partial(x)$ $(c_0, c_1 \in \cc\ \textrm{and}\ x\in L)$. 

\subsection{Example 2}
Let $c\neq 0$ be a complex number and consider functional equation
\begin{equation}
\label{exm03}
(t_1^3+t_2^3)f(x)-(t_1^2+t_2^2)f((t_1+t_2)x)+t_2f(t_1^2x)+t_1f(t_2^2x)=c\cdot x,
\end{equation}
i.e. $n=4$ and $a_1=t_1^3+t_2^3$, $a_2=-(t_1^2+t_2^2)$, $a_3=t_2$, $a_4=t_1$, $\al_1=1$, $\al_2=t_1+t_2$, $\al_3=t_1^2$, $\al_4=t_2^2$, where $t_1$ and $t_2$ are algebraically independent numbers over the rationals, $L=\qq(t_1, t_2)$ and $L\subset K$ such that $K$ is finitely generated and it is contained in the algebraic closure of $L$. First of all compute the possible coefficients of the corresponding linear system of equations:
\begin{equation}
\label{coeffexm02}
\begin{split}
&\sum_{i=1}^4a_i\al_i=0, \\
&\sum_{i=1}^4a_i\partial_1(\al_i)=\sum_{i=1}^4a_i\partial_2(\al_i)=-(t_1-t_2)^2,\\
&\sum_{i=1}^4a_i \partial_1 \partial_2(\al_i)=0, \ \ \sum_{i=1}^4a_i\partial^2_1(\al_i)=2t_2, \ \ \sum_{i=1}^4a_i\partial^2_2(\al_i)=2t_1;\\
\end{split}
\end{equation}
note that any higher order derivative must be zero because $\al_i$'s are at most second order polynomials of the variables $t_1$ and $t_2$. We are going to use Corollary \ref{simpleextension} to compute differential operator solutions
$$D=\sum_{j_1=0}^{J_1}\sum_{j_2=0}^{J_2} c_{j_1 j_2}\partial_1^{j_1} \partial_{2}^{j_2}$$
of degree one, two and three.
\begin{itemize}
\item If $J_1=J_2=1$ then we have four equations for the coefficients of $D$.
\end{itemize} 
Using (\ref{coeffexm02})
\begin{equation*}
\begin{split}
&\sum_{i=1}^4a_iD_{11}(\al_i)=0,\\
&\sum_{i=1}^4a_iD_{10}(\al_i)=-c_{11}(t_1-t_2)^2=0\ \ \Rightarrow\ \ c_{11}=0,\\
&\sum_{i=1}^4a_iD_{01}(\al_i)=0,\\
&\sum_{i=1}^4a_iD_{00}(\al_i)=-(c_{10}+c_{01})(t_1-t_2)^2\neq 0\\
\end{split}
\end{equation*}
and a non-zero particular solution is $(c/\tilde{c})\cdot D$, where
$$D=c_{00}+c_{10}\partial_1+c_{01}\partial_2,$$
$c_{00}, c_{10}$ and $c_{01}\in \cc$ such that $\tilde{c}:=-(c_{10}+c_{01})(t_1-t_2)^2\neq 0$.
\begin{itemize}
\item If $J_1=J_2=2$ then we have nine equations for the coefficients of $D$. 
\end{itemize} 
Using (\ref{coeffexm02})
\begin{equation*}
\begin{split}
&\sum_{i=1}^4a_iD_{22}(\al_i)=0,\\
&\sum_{i=1}^4a_iD_{21}(\al_i)=-2c_{22}(t_1-t_2)^2=0\ \ \Rightarrow\ \ c_{22}=0\\
&\sum_{i=1}^4a_iD_{20}(\al_i)=-c_{21}(t_1-t_2)^2=0\ \ \Rightarrow\ \ c_{21}=0,\\
&\sum_{i=1}^4a_iD_{12}(\al_i)=0,\\
&\sum_{i=1}^4a_iD_{11}(\al_i)=-2c_{12}(t_1-t_2)^2=0\ \ \Rightarrow\ \ c_{12}=0,\\
&\sum_{i=1}^4a_iD_{10}(\al_i)=-\left(2c_{20}+c_{11}\right)(t_1-t_2)^2=0\ \ \Rightarrow\ \ c_{11}=-2c_{20},\\
&\sum_{i=1}^4a_iD_{02}(\al_i)=0,\\
\end{split}
\end{equation*}
\begin{equation*}
\begin{split}
&\sum_{i=1}^4a_iD_{01}(\al_i)=-\left(2c_{02}+c_{11}\right)(t_1-t_2)^2=0\ \ \Rightarrow\ \ c_{11}=-2c_{02},\\
&\sum_{i=1}^4a_iD_{00}(\al_i)=-(c_{10}+c_{01})(t_1-t_2)^2-c_{11}(t_1+t_2)\neq 0
\end{split}
\end{equation*}
and a non-zero particular solution is $(c/\tilde{c})\cdot D$, where
$$D=c_{00}+c_{10}\partial_1+c_{01}\partial_2-\frac{c_{11}}{2}\left(\partial_1^2-2\partial_1 \partial_2 +\partial_2^2\right),$$
$c_{00}, c_{10}, c_{01}$ and $c_{11}\in \cc$ such that 
$$\tilde{c}:=-(c_{10}+c_{01})(t_1-t_2)^2-c_{11}(t_1+t_2)\neq 0.$$
\begin{itemize}
\item If $J_1=J_2=3$ then we have sixteen equations for the coefficients of $D$. 
\end{itemize} 
Using (\ref{coeffexm02})
\begin{equation*}
\begin{split}
&\sum_{i=1}^4a_iD_{33}(\al_i)=0,\\
&\sum_{i=1}^4a_iD_{32}(\al_i)=-3c_{33}(t_1-t_2)^2=0\ \ \Rightarrow\ \ c_{33}=0\\
&\sum_{i=1}^4a_iD_{31}(\al_i)=-2c_{32}(t_1-t_2)^2=0\ \ \Rightarrow\ \ c_{32}=0,\\
&\sum_{i=1}^4a_iD_{30}(\al_i)=-c_{31}(t_1-t_2)^2=0\ \ \Rightarrow\ \ c_{31}=0,\\
&\sum_{i=1}^4a_iD_{23}(\al_i)=0,\\
&\sum_{i=1}^4a_iD_{22}(\al_i)=-3c_{23}(t_1-t_2)^2=0\ \ \Rightarrow\ \ c_{23}=0,\\
&\sum_{i=1}^4a_iD_{21}(\al_i)=-2c_{22}(t_1-t_2)^2=0\ \ \Rightarrow\ \ c_{22}=0,\\
&\sum_{i=1}^4a_iD_{20}(\al_i)=-(c_{21}+3c_{30}(t_1-t_2)^2=0\ \ \Rightarrow\ \ c_{21}=-3c_{30},\\
&\sum_{i=1}^4a_iD_{13}(\al_i)=0,\\
&\sum_{i=1}^4a_iD_{12}(\al_i)=-3c_{13}(t_1-t_2)^2=0\ \ \Rightarrow\ \ c_{13}=0\\
&\sum_{i=1}^4a_iD_{11}(\al_i)=-2(c_{12}+c_{21})(t_1-t_2)^2=0\ \ \Rightarrow\ \ c_{12}=-c_{21},\\
&\sum_{i=1}^4a_iD_{10}(\al_i)=-(c_{11}+2c_{20})(t_1-t_2)^2+2t_1c_{12}+6t_2c_{30}=0\ \ \Rightarrow\ \ \\
&(c_{11}+2c_{20})(t_1-t_2)^2=2t_1c_{12}+6t_2c_{30},\\
\end{split}
\end{equation*}
\begin{equation*}
\begin{split}
&\sum_{i=1}^4a_iD_{03}(\al_i)=0,\\
&\sum_{i=1}^4a_iD_{02}(\al_i)=-(c_{12}+3c_{03})(t_1-t_2)^2=0\ \ \Rightarrow\ \ c_{12}=-3c_{03},\\
&\sum_{i=1}^4a_iD_{01}(\al_i)=-(c_{11}+2c_{02})(t_1-t_2)^2+2t_2c_{21}+6t_1c_{03}=0\ \ \Rightarrow\ \ \\
&(c_{11}+2c_{02})(t_1-t_2)^2=2t_2c_{21}+6t_1c_{03},\\
&\sum_{i=1}^4a_iD_{00}(\al_i)=-(c_{01}+c_{10})(t_1-t_2)^2+2t_1c_{02}+2t_2c_{20}\neq 0\\
\end{split}
\end{equation*}
and a non-zero particular solution is $(c/\tilde{c})\cdot D$, where
$$D=c_{00}+c_{10}\partial_1+c_{01}\partial_2-\frac{c_{11}}{2}\left(\partial_1^2-2\partial_1 \partial_2 +\partial_2^2\right)+\frac{t_1+t_2}{(t_1-t_2)^2}c_{12}\left(\partial_1^2-\partial_2^2\right)+$$
$$\frac{c_{12}}{3}\left(\partial_1^3+3\partial_1 \partial_2^2-3\partial_2\partial_1^2-\partial_2^3\right),$$
$c_{00}, c_{10}, c_{01}, c_{11}$ and $c_{12}\in \cc$ such that 
$$\tilde{c}:=-(c_{10}+c_{01})(t_1-t_2)^2-c_{11}(t_1+t_2)-2\frac{t_1+t_2}{t_1-t_2}c_{12}\neq 0.$$
The space of the solutions on $K$ is $D+S_1^0.$ Equation $\displaystyle{\sum_{i=1}^4 a_i\al_i}=0$ implies that any exponential function in $S_1$ belongs to $S_1^0$, i.e. 
\begin{equation}
\label{huge}
(t_1^3+t_2^3)-(t_1^2+t_2^2)(s_1+s_2)+t_2\cdot s_1^2+t_1\cdot s_2^2=0,
\end{equation}
where $\phi(t_1)=s_1$ and $\phi(t_2)=s_2$ are algebraically independent over the rationals. Since the coefficient $t_2/t_1$ of the normalized characteristic polynomial
$$p(x,y)=\frac{t_1^3+t_2^3}{t_1}-\frac{t_1^2+t_2^2}{t_1}(x+y)+\frac{t_2}{t_1}\cdot x^2+y^2$$
is transcendent it has algebraically independent roots; see \cite{V2} and \cite{VV}. Especially $\phi(t_1)=t_1$ and $\phi(t_2)=t_2$ are solutions of equation \eqref{huge}. To find the generating elements of the form $\phi\circ D$ in $S_1^0$ we need to use
\eqref{diffophom}:
$$\sum_{i=1}^4 \phi^{-1}(a_i)D_{m_1 m_2}(\al_i)=0\ \ (0\leq m_1 \leq J_1, 0\leq m_2 \leq J_2).$$
If $\phi(t_1)=t_1$ and $\phi(t_2)=t_2$ then $\phi^{-1}(a_i)=a_i$ ($i=1, \ldots, 4$). Therefore $\phi\circ D\in S_1^0$, where $D$ is one of the differential operators above provided that the  surviving coefficients is choosen such that $\tilde{c}=0$. For example, a differential operator of degree one in $S_1^0$ is
$$D=c_{00}+c_{10}\partial_1+c_{01}\partial_2,\ \ \textrm{where}\ \ c_{10}+c_{01}=0.$$

\section{Spectral synthesis in the variety containing monomial solutions of higher degree}

In what follows we are going to give a survey of the higher order version of section 3. By Theorem 6.3 in \cite{KL} the space $S_p$ is spanned by the functions of the form
$$(x_1, \ldots, x_p)\mapsto \phi_1\circ D_1(x_1)\cdot \ldots \cdot \phi_p \circ D_p(x_p),$$
where $\phi:=\phi_1 \cdot \ldots \cdot \phi_p$ is an exponential element in $S_p$ and $D_1$, $\ldots$, $D_p$ are differential operators on $K$, where $K$ is a finitely generated field containing the parameters $\al_i$ and $\be_i$ $(i=1, \ldots, n)$. 

\subsection{The case of transcendence degree $0$ (algebraic parameters)}

If the parameters $\alpha_1$, $\ldots$, $\alpha_n$, $\be_1$, $\ldots$, $\be_n$ are algebraic numbers (the case of transcendence degree $0$) then we have no non-trivial differential operators on $K=\qq(\al_1, \ldots, \al_n,\be_1, \ldots, \be_n)$; cf. subsection 3.2. The analogue results can be easily formulated for $p>1$ by using the higher order version of Theorem \ref{t1o1}; see also section 4 in \cite{KV2}. 

\subsection{The case of higher transcendence degree}

Suppose that the transcendence degree of the parameters $\alpha_1$, $\ldots$, $\alpha_n$, $\be_1$, $\ldots$, $\be_n$ is $k$, i.e. we have a field extension $L=\mathbb{Q}(t_1, \ldots, t_k)$ such that $t_1$, $\ldots$, $t_k$ are algebraically independent over the rationals, $L\subset K$ and $K$ is contained in the algebraic closure of $L$; especially it contains the parameters $\alpha_i$'s and $\beta_i$'s. 

\begin{lemma}
\label{higher1} Suppose that $\phi(x)=x$ $(x\in K)$; then
$$\phi\circ D(\al_ix)=\sum_{m_1=0}^{J_1}\ldots \sum_{m_k=0}^{J_k}D^{\phi}_{m_1 \ldots m_k}(\al_i)\partial_1^{m_1}\ldots \partial_k^{m_k}(x)$$
for any $x\in K$. 
\end{lemma}

The proof is a straightforward calculation by formula \eqref{eval1}. The key step of the application of the higher order spectral synthesis is to formulate the necessary and sufficient conditions for 
$$\tilde{F}_p(x_1, \ldots, x_p):=\phi_1\circ D_1(x_1)\cdot \ldots \cdot \phi_p \circ D_p(x_p)$$
to satisfy system (\ref{multimultitilde}). In what follows we discuss only the case of $p=2$. The case of $p>2$ can be investigated in a similar way because of the inductive argument as follows. The first equation of system  (\ref{multimultitilde}) implies that
\begin{equation}\label{higher2}
\sum_{i=1}^na_i \phi_1\circ D_1(\al_i x_1)\cdot \phi_2 \circ D_2(\al_i x_2)=\tilde{c}\cdot x_1\cdot x_2
\end{equation}
In order to use the results in the previous sections let the variables $x_2$ be considered as a non-zero given constant in $K$; the most simple choice is $x_2=1$. Then 
\begin{equation}\label{higher3}
\sum_{i=1}^n\underbrace{a_i \phi_2\circ D_2(\al_i x_2)}_{\textrm{new coefficients}}\cdot \phi_1\circ D_1(\al_i x_1)=\tilde{c} \cdot x_1 \cdot x_2,
\end{equation}
where the new constant is $\tilde{c}\cdot x_2$, i.e. $\phi_1\circ D_1$ satisfies one of the conditions to be the additive solution of an inhomogeneous equation with some new coefficients. If $\tilde{c}\neq 0$ then we have, by Proposition \ref{ldiffuj}, that $\phi_1(x)=x$ ($x\in K$) and
\begin{equation}\label{higher4}
\begin{split}
&\sum_{i=1}^na_i \phi_2\circ D_2(\al_i x_2)\cdot \left(D_{1}^{\phi_1}\right)_{m_{11}\ldots m_{1k}}(\al_i)=0\ \ \textrm{if} \ \ m_{11}+\ldots+m_{12}\geq 1\\
&\sum_{i=1}^na_i \phi_2\circ D_2(\al_i x_2)\cdot D_{1}^{\phi_1}(\al_i)=\tilde{c}\cdot x_2.
\end{split}
\end{equation}
In a similar way, $\phi_2(x)=x$ ($x\in K$) and, by using Lemma \ref{higher1}, the left hand side of the first equation is the action of a differential operator at $x_2$ for any given $m_{11}, \ldots, m_{1k}$:
{\footnotesize{$$\sum_{m_{21}=0}^{J_{21}}\ldots \sum_{m_{2k}=0}^{J_{2k}}\sum_{i=1}^na_i \left(D_2^{\phi_2}\right)_{m_{21} \ldots m_{2k}}(\al_i)\cdot \left(D_{1}^{\phi_1}\right)_{m_{11}\ldots m_{1k}}(\al_i)\cdot \partial_1^{m_{21}}\ldots \partial_k^{m_{2k}}(x)=0.$$}}

\noindent
Therefore
\begin{equation}\label{higher5}
\begin{split}
&\sum_{i=1}^na_i \left(D_{1}^{\phi_1}\right)_{m_{11}\ldots m_{1k}}(\al_i)\cdot \left(D_{2}^{\phi_2}\right)_{m_{21}\ldots m_{2k}}(\al_i)=0\\
& \textrm{if} \ \ m_{11}+\ldots+m_{1k}\geq 1\ \ \textrm{and}\ \ 0\leq m_{21}\leq J_{21}, \ldots, 0\leq m_{2k}\leq J_{2k}.\\
\end{split}
\end{equation}
On the other hand the second equation of \eqref{higher4} shows that $\phi_2\circ D_2$ satisfies one of the conditions to be the additive solution of an inhomogeneous equation with some new coefficients. By Proposition \ref{ldiffuj} 
\begin{equation}\label{higher6}
\begin{split}
&\sum_{i=1}^na_i \left(D_{2}^{\phi_2}\right)_{m_{21}\ldots m_{2k}}(\al_i)\cdot D_{1}^{\phi_1}(\al_i)=0\ \ \textrm{if} \ \ m_{21}+\ldots +m_{2k}\geq 1\\
&\sum_{i=1}^na_i D_2^{\phi_2}(\al_i)\cdot D_{1}^{\phi_1}(\al_i)=\tilde{c}.
\end{split}
\end{equation}
Equations (\ref{higher5}) and (\ref{higher6}) give that
\begin{equation}\label{higher7}
\begin{split}
&\sum_{i=1}^na_i \left(D_{1}^{\phi_1}\right)_{m_{11}\ldots m_{1k}}(\al_i)\cdot \left(D_{2}^{\phi_2}\right)_{m_{21}\ldots m_{2k}}(\al_i)=0\\
&\textrm{if} \ \ m_{11}+\ldots+m_{1k}+m_{21}+\ldots +m_{2k}\geq 1\\
&\sum_{i=1}^na_i D_1^{\phi_1}(\al_i)\cdot D_{2}^{\phi_2}(\al_i)=\tilde{c}.\\
\end{split}
\end{equation}
System \eqref{multimultitilde} implies similar equations containing the parameters $\beta_1$, $\ldots$, $\beta_n$ (pure case) or $\al_1$, $\ldots$, $\al_n$, $\be_{1}$, $\ldots$, $\be_{n}$ (mixed cases). According to the definition of $D_{m_1 \ldots m_k}$ the result is a non-linear (polynomial; especially quadratic because of $p=2$) system \eqref{higher7} for the quantities $c'_{j_{11} \ldots j_{1k}}$'s and $c'_{j_{21} \ldots j_{2k}}$'s, where
$$D_1^{\phi_1}=\sum_{j_{11}=0}^{J_{11}}\ldots \sum_{j_{1k}=0}^{J_{1k}}c_{j_{11} \ldots j_{1k}}' \partial_1^{j_{11}}\ldots \partial_k^{j_{1k}},$$
$$D_2^{\phi_2}=\sum_{j_{21}=0}^{J_{21}}\ldots \sum_{j_{2k}=0}^{J_{2k}}c_{j_{21} \ldots j_{2k}}' \partial_1^{j_{21}}\ldots \partial_k^{j_{2k}}, $$
$$c'_{j_{11} \ldots j_{1k}}=\phi_1(c_{j_{11} \ldots  j_{1k}}), \ \ c'_{j_{21} \ldots j_{2k}}=\phi_2(c_{j_{21} \ldots j_{2k}}).$$

\subsubsection{The case of transcendence degree $1$}

For the sake of simplicity consider the case of $k=1$ (transcendence degree one), $p=2$ and the separated version of the functional equations as in section 4 (see also Remark \ref{separation} and Remark \ref{separation1}). For any fixed index $m_{21}$ system (\ref{higher7}) results in a diagonal form like (\ref{diagonal1}) with some new coefficients:
\begin{equation}\label{higher8}
\begin{split}
&\sum_{i=1}^n\underbrace{a_i \cdot \left(D_{2}^{\phi_2}\right)_{m_{21}}(\al_i)}_{\textrm{new coefficients}}\cdot \left(D_{1}^{\phi_1}\right)_{m_{11}}(\al_i)=0\ \ \textrm{if} \ \ m_{11}+m_{21}\geq 1\\
&\sum_{i=1}^n\underbrace{a_i \cdot D_{2}^{\phi_2}(\al_i)}_{\textrm{new coefficients}}\cdot D_1^{\phi_1}(\al_i)=\tilde{c}.\\
\end{split}
\end{equation}
Therefore we can conclude that
\begin{equation}\label{higher9}
\begin{split}
&\sum_{i=1}^na_i \al_i \left(D_{2}^{\phi_2}\right)_{m_{21}}(\al_i)=0\\
&\sum_{i=1}^na_i \partial(\al_i) \left(D_{2}^{\phi_2}\right)_{m_{21}}(\al_i)=0\\
& \ldots\\
&\sum_{i=1}^na_i \partial^{J_{11}-1}(\al_i) \left(D_{2}^{\phi_2}\right)_{m_{21}}(\al_i)=0\\
&\sum_{i=1}^na_i \partial^{J_{11}}(\al_i) \left(D_{2}^{\phi_2}\right)_{m_{21}}(\al_i)=\left\{\begin{array}{rl}
&0\ \ \textrm{if} \ \ m_{21}\geq 1\\
&\\
&\tilde{c}\ \ \textrm{if}\ \ m_{21}=0
\end{array}
\right.\\
\end{split}
\end{equation}
Taking the terms of the form $a_i\partial^{j_{11}}(\al_i)$ ($i=1, \ldots, n$) as new coefficients in the $j_{11}$th equation we can give similar conclusions by (\ref{diagonal1}) in case of each equation above:
\begin{equation}\label{higher10}
\begin{split}
&\sum_{i=1}^na_i \al_i^2=0\ \ (j_{11}+j_{21}=0),\\
&\ldots \\
&\sum_{i=1}^na_i \partial^{j_{11}}(\al_i) \partial^{j_{21}}(\al_i)=0\ \ (j_{11}+j_{21}< J_{11}+J_{21}),\\
&c_{J_{11}}'c_{J_{21}}'\sum_{i=1}^na_i \partial^{J_{11}}(\al_i) \partial^{J_{21}}(\al_i)=\tilde{c} \ \ (j_{11}+j_{21}= J_{11}+J_{21}).\\
\end{split}
\end{equation}
These are necessary and sufficient conditions in case of transendence degree $1$ for the generating elements of the space $S_2$. The coefficients can be arbitrarily choosen except the greatest one in both $D_1$ and $D_2$ provided that system (\ref{higher10}) holds.

\subsubsection{An explicite example}  Let $c\neq 0$ be a complex number and consider functional equation 
\begin{equation}
\label{exm4}
t^6f(tx)-t^4f(t^2x)-t^2f(t^3x)+f(t^4x)=c\cdot x^2 \ \ 
\end{equation}
i.e. $n=4$ and $a_1=t^6$, $a_2=-t^4$, $a_3=-t^2$, $a_4=1$, $\al_1=t$, $\al_2=t^2$, $\al_3=t^3$, $\al_4=t^4$, where $t$ is a transcendental number over the rationals, $L=\qq(t)$ and $L\subset K$ such that $K$ is a finitely generated field and it is contained in the algebraic closure of $L$. Using that
\begin{equation}
\label{exm4no1}
\begin{split}
&\sum_{i=1}^4 a_i\al_i^2=t^6\cdot t^2-t^4\cdot t^4-t^2\cdot t^6+t^8=0\\
&\sum_{i=1}^4 a_i \al_i \partial (\al_i)=t^7-t^4\cdot t^2\cdot(2t)-t^2\cdot t^3\cdot (3t^2)+t^4\cdot (4t^3)=0\\
\end{split}
\end{equation}
but 
\begin{equation}
\label{exm4no2}
\begin{split}
&\sum_{i=1}^4 a_i \partial(a_i)\partial(\al_i)=t^6- t^4\cdot(2t)^2-t^2\cdot (3t^2)^2+(4t^3)^2=4t^6 \neq 0\\
\end{split}
\end{equation}
it follows that the generating elements of $S_2$ are the products of first order exponential monomials $D_1^{\phi_1}(x)=c_{10}'x+c_{11}'\partial(x)\ \ \textrm{and}\ \ D_2^{\phi_2}(x)=c_{20}'x+c_{21}'\partial(x)$
such that $c_{11}'\cdot c_{21}'=\frac{c}{\ 4t^6}$, where  $c_{11}'=\phi_1(c_{11})$, $c_{21}'=\phi_2(c_{21})$  and the differential operators are of the form $D_1(x)=c_{10}x+c_{11}\partial(x)\ \ \textrm{and}\ \ D_2(x)=c_{20}x+c_{21}\partial(x).$ 

\section{Concluding remarks} The spectral synthesis in $S_p^*$ allows us to describe the entire space of solutions of an inhomogeneous linear functional equation on a large class of finitely generated fields, at least theoretically. In general we need to solve inhomogeneous linear systems of equations to find the solutions. The conclusions for the homogeneous case has been also formulated step by step in terms of some remarks. The discussion of the case of algebraic parameters (the transcendence degree is zero) is relatively simple because of the trivial action of any differential operator on algebraic numbers. The case of transcendence degree one gives the best results independently of the degree of the monomial solutions. In case of higher degree of transcendence the problem becomes much more difficult because of the missing main terms of the differential operators and the increasing number of the equations. Some explicite examples are also presented to illustrate both the effectivity and the difficulties of the method in practice. Computer assisted methods can be successful for the formulation and the solution of large systems of linear equations belonging to the case of higher order monomial solutions and/or higher degree of transcendence. 

\end{document}